\documentclass[12pt]{article}

\usepackage[T2A]{fontenc}

\usepackage{amsfonts}
\usepackage{amssymb}
\usepackage{amsmath}
\usepackage{amsthm}

\def \le {\leqslant}
\def \ge {\geqslant}

\begin{document}
\begin{Large}
\centerline{\bf A note on badly approximable affine forms and winning sets} \centerline{ \bf N.G. Moshchevitin\footnote{ Research is supported by RFBF grant No 09-01-00371a
 }  
  }
\end{Large}

 \vskip+0.7cm
 \centerline{Abstract}

\begin{small}
 We prove a result on  inhomogeneous Diophantine approximations related to  the theory of $(\alpha,\beta)$-games.
\end{small}
 \vskip+0.7cm

{ \bf 1. Introduction.}\,\,\, J.W.S. Cassels in his book \cite{Cassil}  (Theorem 10, Chapter 5)  describes the following result  on inhomogeneous linear Diophantine approximations.

{\bf Theorem A. }\,\,\, {\it Given positive integers $n,m$ there exists a positive constant
$\Gamma_{m,n} $ with the following property. Let
$$\theta_{i,j},\,\,\, 1\le i \le m,\,\,\, 1\le j\le n$$ be real numbers and let
$$
L_j({\bf x}) = \sum_{i=1}^m \theta_{i,j}x_i ,\,\,\ 1\le j\le n
$$ be the corresponding system of linear homogeneous forms in integer variables $x_1,...,x_m$. Then there exist a collection of real numbers $\eta_1,...,\eta_n$ such that
 $$\inf_{{\bf x} \in \mathbb{Z}^m\setminus\{{\bf 0}\}}
\left( \max_{1\le j\le n}||L_j({\bf x}) -\eta_j ||\right)^n \left(\max_{1\le i \le m} |x_i|\right)^m >\Gamma_{m,n}.
$$}
Here and in the sequel $||\cdot ||$ denotes the distance to the nearest integer.

The method of  the proof of this result
generalizes a
 construction introduced by A.Khintchine in \cite{H}, \cite{a} (see also \cite{HH}).

Theorem A is a particular case of a more general result proved by V.Jarnik in 1941 (see \cite{J41},\cite{J41bis}).

Given a collection of real numbers
$$
\Theta =\{
\theta_{i,j},\,\,\, 1\le i \le m,\,\,\, 1\le j\le n
\}$$
consider  the
 function
$$
\psi_{\Theta} (t)  =
 \,\,\,\,\,\min_{(y_1,...,y_n)\in\mathbb{R}^n: \, 0<\max_{1\le j\le n}|y_j| \le t
}\,\,\,\,\,\max_{1\le i\le m}|| \sum_{1\le j \le n}\theta_{i,j}y_j||
.
$$

{\bf Theorem B. }(V. Jarnik \cite{J41},\cite{J41bis})\,\,\, {\it Given positive integers $n,m$ there exists a positive constant
$\Gamma_{m,n}^* $ with the following property. 
Suppose $\psi (t) $
to be a function decreasing to zero as $t\to+\infty$. Let   $\rho (t)$ 
be the function inverse to the function
 $t\mapsto 1/\psi (t)$.
 Let for all  $t$ large enough one has
 $$
\psi_{\Theta}\le \psi (t).$$
Then there exists a vector 
$$
(\eta_1,...,\eta_n) $$
such that
$$\inf_{{\bf x} \in \mathbb{Z}^m\setminus\{{\bf 0}\}}  
\left( \max_{1\le j\le n}||L_j({\bf x}) -\eta_j ||\right)\cdot \rho \left(\max_{1\le i\le m}||x_i||\right) >\Gamma_{m,n}^*
.$$
 }

To obtain Theorem A one must take $\psi (t) = t^{-n/m}$  in Theorem B.

Theorem A was generalized by D. Kleinbock in \cite{KL}.

{\bf Theorem C. }(D. Kleinbock,\cite{KL})\,\,{\it Define ${\cal B}$ to be a set of real $(m+1)\times n$ martices
$$(\theta_{i,j},\eta_j ),\,\,\, 1\le i \le m,\,\,\, 1\le j\le n;\,\,\,\,
$$
for which    the inequality
$$\inf_{{\bf x} \in \mathbb{Z}^m\setminus\{{\bf 0}\}}
\left( \max_{1\le j\le n}||L_j({\bf x}) -\eta_j ||\right)^n \left(\max_{1\le i \le m} |x_i|\right)^m >0.
$$ is valid. Then ${\cal B}$
has full Hausdorff dimension in $\mathbb{R}^{mn+n}$.}

 Kleinbock's proof was based on the consideration of special flows on homogeneous spaces.
Y. Bugeaud, S. Harrap, S. Kristensen and S. Velani \cite{BU} give a simple proof of Kleinbock's theorem. Moreover they establish the following result.

{\bf Theorem D. }(see \cite{BU})\,\,\, {\it For any collection  of real numbers
$$
\Theta =\{
\theta_{i,j},\,\,\, 1\le i \le m,\,\,\, 1\le j\le n
\}$$
consider the set ${\cal B}(\Theta )$ of real vectors
$$
(\eta_1,...,\eta_n) $$
for which
$$\inf_{{\bf x} \in \mathbb{Z}^m\setminus\{{\bf 0}\}}
\left( \max_{1\le j\le n}||L_j({\bf x}) -\eta_j ||\right)^n \left(\max_{1\le i \le m} |x_i|\right)^m >0.
$$
Then ${\cal B}(\Theta)$
has full Hausdorff dimension in $\mathbb{R}^{n}$.}

A remarkable improvement was obtained recently by J. Tseng.

{\bf Theorem  E }(J.Tseng, \cite{Ts})\,\,\, {\it For any $\theta$ the set  ${\cal B}$ of real
numbers $\eta$ such that
$$
\inf_{x\in \mathbb{Z}\setminus\{ 0\}} |x|\cdot ||\theta x+\eta||>0
$$
is an $\alpha$-winning set  for any $\alpha \in (0,1/8)$.}

Moreover as it was mentioned in Remark 2.3 in \cite{Ts}  J. Tseng and M. Einsiedler obtained a generalization of Theorem E to
the case of arbitrary systems of linear forms.

We do not want to discuss the definition of winning set and the metrical properties of such sets. We refer to the book  \cite{SCH1} and  the
paper \cite{SCH} by W.M. Schmidt where all necessary definitions and results are given. We only want to note that any winning set in
$\mathbb{R}^n$ has full Hausdorff dimension. So Theorem D follows from the result by  J. Tseng and M.Einsiedler.

 Also we would like to say that there exist various papers devoted to the study of winning sets and their properties (see
for example \cite{fb},\cite{fb1},\cite{ah},\cite{ah1},\cite{dre},\cite{KLEI},\cite{Mos},\cite{SCH1},\cite{SCH},\cite{SCH2}).

In this note we would like to announce an improvement of Theorem D.

{\bf Theorem 1.} \,\,{\it Let $\alpha \in (0,1/2)$. Then for any collection  of real numbers
$$
\Theta =\{
\theta_{i,j},\,\,\, 1\le i \le m,\,\,\, 1\le j\le n
\}$$
 the set ${\cal B}(\Theta )$ of real vectors
$$
(\eta_1,...,\eta_n) $$
for which
$$\inf_{{\bf x} \in \mathbb{Z}^m\setminus\{{\bf 0}\}}
\left( \max_{1\le j\le n}||L_j({\bf x}) -\eta_j ||\right)^n\left(\max_{1\le i \le m} |x_i|\right)^m >0
$$
is an $\alpha$-winning set in $\mathbb{R}^n$.}

The followong result generalizes Theorem 1. 

{\bf Theorem 2.} \,\,{\it Let  $\alpha \in (0,1/2)$. 
Suppose $\psi (t) $
to be a function decreasing to zero as $t\to+\infty$. Let   $\rho (t)$ 
be the function inverse to the function
 $t\mapsto 1/\psi (t)$.
 Let for all  $t$ large enough one has
 $$
\psi_{\Theta}\le \psi (t).$$
Then the set 
${\cal B}(\Theta )$ of all vectors 
$$
(\eta_1,...,\eta_n) $$
such that
$$\inf_{{\bf x} \in \mathbb{Z}^m\setminus\{{\bf 0}\}}  
\left( \max_{1\le j\le n}||L_j({\bf x}) -\eta_j ||\right)\cdot \rho \left(\max_{1\le i\le m}||x_i||\right) >0
$$
is an   $\alpha$-winning set in  $\mathbb{R}^n$.

 }

On cas easily see that 
Theorem 1 is a particular case of Theorem 2.
In the next section we sketch the main steps of the proof of Theorem 1.
The proof of Theorem 2 follows the steps of the proof of Theorem 1.

{\bf 2. Sketch of the proof for Theorem 1.}

{\bf Lemma 1.}\,\,\, {\it Let a sequence of reals $t_r, r =1,2,3,...$ satisfy  the lacunarity condition
\begin{equation}
\frac{t_{r+1}}{t_r} \ge M
 ,\,\,\,\,\,
r=1,2,3,...
 \label{lac}
\end{equation}
for some $ M>1$. Let a sequence $\Lambda\subset \mathbb{Z}^n$ of integer vectors ${\bf u}^{(r)}= (u^{(r)}_1, ....,u^{(r)}_n  ) \in \mathbb{Z}^n$
be such that
\begin{equation}
t_r^2= (u^{(r)}_1)^2+ ....+(u^{(r)}_n)^2. \label{vec}
\end{equation}
Then the set
$$
N(\Lambda ) =\{ \eta =(\eta_1,...,\eta_n)\in\mathbb{R}^n:\,\,\, \exists \, c(\eta)>0 \,\text{such that}\, || u^{(r)}_1\eta_1+ ...+u^{(r)}_n
\eta_n||\ge c(\eta )\,\,\ \forall r \in \mathbb{N} \}
$$
is $\alpha$-winning for any $\alpha \in (0,1/2)$.}

Lemma 1 is a generalization of Lemma 2 from \S 6, Chapter 5 of Cassels' book \cite{Cassil}.  In order to prove Lemma 1 one should follow the
proof of Theorem from the paper \cite{Mos}. This proof is a generalization of Theorem 4 proof from \cite{SCH}. In this proof two main arguments
are used: the original Schmidt's lemma on escaping from dangerous points (Lemma 15 from \cite{SCH}) and dichotomy process.  To construct a proof
of Lemma 1 from the proof of Theorem from \cite{Mos} in the case $n=1$ one needs only  to make minor corrections in the choise of the parameters
(such as put $\varepsilon = 0$). In the multidimensional case one must use Schmidt's multidimensional escaping lemma (Lemma 1B, Chapter 3 form
\cite{SCH1}) instead of the one-dimensional escaping lemma (Lemma 15 from \cite{SCH}). This does not enable to make dichotomy process in the
proof. So we give a complete proof of Lemma 1 in the next section.

We would like to note that Cassels' lemma was recently improved by I. Rochev \cite{RO} by means of Peres-Shlag's method \cite{PS}.

To deduce Theorems 1,2 from Lemma 1 we must follow the arguments of the proof of Theorem 10 from \S 6  Chapter 5 of the book \cite{Cassil}.
Given numbers $\theta_{i,j},\,\,\, 1\le i \le m,\,\,\, 1\le j\le n$ we take the sequence of vectors ${\bf u}^{(r)}$ form Lemma 4 \S 6 Chapter 5
from \cite{Cassil}
 satisfying the
conditions ({\ref{lac},\ref{vec}) with $M=3$. Then we take arbitrary $\eta$ from the set $N(\Lambda )$. Following the proof of Theorem 10 from \S
6 Chapter 5 of the book \cite{Cassil} we see that for any integer vector ${\bf x}=(x_1,...,x_m) \in \mathbb{Z}^m\setminus \{{\bf 0}\}$ one has
$$
\left( \max_{1\le j\le n}||L_j({\bf x}) -\eta_j ||\right)^n \left(\max_{1\le i \le m} |x_i|\right)^m >0.
$$

{\bf 3. Proof of Lemma 1.}\,\,\,

As usual for $\alpha, \beta \in (0,1)$ put $\gamma = 1+\alpha\beta -2\alpha >0$.

For a ball $B\subset \mathbb{R}^n$ with the center $O$  and radius $\rho$ we consider its boundary $S = \partial B$. Let $\mu $ be normalized
Lebesgue measure on $S$. So $ \int_S d\mu = \mu S = 1$. Let $ x \in S$. Define $\pi (x)\subset \mathbb{R}^n$ to be the $n-1$-dimensional affine
subspace passing through the center $O$ of the ball $B$ orthogonal to one-dimensional subspace passing through the point $O$ and $x$. Define
$\Pi (x)$ to be the "half-space" with $\pi (x)$ as a boundary such that $x \in \Pi (x)$.

 Given $\alpha, \beta \in (0,1)$ consider the unique "halfspace" $ \Pi_{\alpha, \beta,\rho }(x)$ such that
$ \Pi_{\alpha, \beta,\rho }(x)\subset \Pi (x)$ and  the distance from  $ \Pi_{\alpha, \beta,\rho }(x)$
 to
 $O$ is equal to $\frac{\gamma \rho}{2}$. Put
 $$\Omega (x) =S\cap \Pi_{\alpha, \beta,\rho }(x),\,\,\,
 \Omega^*(x) = \bigcup_{y\in S:\,\, \Pi (y) \supset \Omega (x)} \{y\}.
 $$
 Obviously the value $\mu \Omega^* (x)$ does not depend on $x\in S$.
 Put
\begin{equation}
 \omega = \omega (\alpha, \beta ) = \mu \Omega^* (x) \in (0,1).
 \label{om}
\end{equation}

 {\bf Lemma 2.}\,\,\,{\it Let $\pi_1, ... , \pi_k$ be any $n-1$-dimensional affine subspaces. Then there exists a point $x\in S$
 such that
 $$\Omega (x) \cap \pi_j = \varnothing,\,\,\,
 $$
for al least $  \lceil\omega k \rceil$ indices $j$.}

Proof. Define $\pi_j'$ to be  a $n-1$-dimensional affine subspace parallel to $\pi_j$ such that $ 0\in \pi_j'$. We take a point $\xi_j\in S$
such that the line passing through  $\xi_j$ and $o$ is orthogonal to $\pi_j'$ and $\xi_j,\pi_j$ lie on the different sides from $\pi_j'$. Note
that if $\xi_j \in \Omega^*(x) $ then $ \pi_j \cap \Omega (x) =\varnothing$. Then define $ \chi_x (y)$ to be the characteristic function of
$\Omega^*(x)$. It should be sufficient to prove that there exists $ x \in S$ such that
$$
f(x) =\sum_{j=1}^k \chi_x(\xi_j) \ge \omega k.
$$
This is true as
$$
\int_S f(x) d\mu (x) = k\int_S \chi_x(\xi) d\mu (x) = k\omega.
$$
Lemma is proved.

 {\bf Lemma 3.}(W.M.Schmidt's  escaping lemma, Lemma 1B, Chapter 3 form \cite{SCH1})  \,\,\,{\it
 Let $t$ be such that
 $$
 (\alpha \beta )^t <\frac{\gamma}{2}.
 $$
 Suppose   a ball $ B_j$ occurs in the game (as a Black ball). Suppose $\pi$ is an $n-1$-dimensional affine subspace passing through the
 center of the ball $B_j$. Then White  can play in such a way that the ball $B_{k+t}$ is contained in the "halfspace" $ \Pi_{\alpha, \beta,\rho_j }(x)$
 such that the boundary of
$ \Pi_{\alpha, \beta,\rho_j }(x)$ is parallel to the subspace $\pi$.
 }

Put
\begin{equation}
t = t(\alpha,\beta, ) =\left\lceil  \frac{\log(\gamma/2)}{\log (\alpha\beta )}\right\rceil,\,\,\,\,\, \tau_k =
 t\cdot\left\lceil \frac{\log k}{\log\left(\frac{1}{1-\omega}\right)} \right\rceil
, \label{tau}
\end{equation}
where $\omega$ is defined in (\ref{om}). From Lemmas 2,3 we obtain the following

{\bf Corollary 1.}\,\,\,{\it Let a ball $B_j$ with the radius $\rho_j$ occurs as a Black ball in the game. Let we hawe a collection of affine
subspaces $\pi_i,\, 1\le i \le k$. Then White can play in such a way that for any point $x \in B_{j+\tau_k}$ the distance between $x$ and any of
subspaces $\pi_i,\, 1\le i \le k$ is greater than $\frac{\rho_{j+\tau_k}\gamma}{2}$.}

Now we are ready to show how to prove Lemma 1.   Take a natural number $k=k (\alpha ,\beta,  M)$ such that
\begin{equation}
\tau_k\times \frac{\log (1/(\alpha\beta ) )}{\log M} +2 < k.
  \label{ka}
\end{equation}
We may suppose without loss of generality that in addition to (\ref{lac}) we have the condition
\begin{equation}
\frac{t_{r+1}}{t_r} \le M^2, \,\,\, r=1,2,3,... .
 \label{ka1}
\end{equation}
Suppose that Black  begin the game with the ball $B_0$ with the radius $\rho$. Define natural numbers $1=r_1\le r_1\le ...\le r_j \le ...$ from
the condition
$$
\frac{1}{2t_{r_{j+1}}} >\rho (\alpha \beta )^{j\tau_k}\ge \frac{1}{2t_{r_{j+1}+1}}.
$$
Then
\begin{equation}
\frac{1}{2t_{r_{j+1}}} >\rho (\alpha \beta )^{j\tau_k}\ge \frac{1}{2M^2t_{r_{j+1}}}
 \label{t}
\end{equation}
(we may assume that $ \frac{1}{t_1} >\rho$).

Define
$$
\varepsilon = \frac{\gamma}{4M^{2+k(\alpha,\beta,\gamma)}}.
$$

 Suppose that we know that for every $x \in B_{j\tau_k}$ the distance between $x$ and every subspace
of the form
\begin{equation}
\{ {\bf y}=(y_1,... , y_n)\in \mathbb{R}^n:\,\, u_1^{(r)}y_1+...+ u_n^{(r)}y_n =a\},\,\,\, a\in \mathbb{Z}
 \label{sp}
\end{equation}
is greater than $\frac{\varepsilon}{t_r}$ for any $r$ from the range $ 1\le r \le r_j$.
Left inequality from  (\ref{t})
shows that any ball  $B_\nu$ 
with $ j\tau_k <\nu \le (j+1)\tau_k$ 
has the radius $ \rho (\alpha\beta)^\nu < 1/(2t_{r_{j+1}})$. So given
$\nu$ from this interval  there exists not more than one "dangerous" affine subspace of the form (\ref{sp}) with $r=\nu$.
 Hence in the range  $ j\tau_k
<\nu\le (j+1)\tau_k $ only $<q_j = r_{j+1} -r_j$ "dangerous" affine subspaces may occur. By the lacunarity condition (\ref{lac}) and the
inequalities (\ref{t}) we have
\begin{equation}
M^{q_j}= M^{r_{j+1}-r_j}\le \frac{t_{r_{j+1}}}{t_{r_j}} \le \frac{M^2}{(\alpha \beta)^{\tau_k}}.
\label{new}
\end{equation}
So
$$
q_j< \frac{\tau_k \log (1/(\alpha\beta))}{\log M} +2 <k
$$
by (\ref{ka}). Now Corollary 1 shows that White can play in such a way that for every $x$ from the ball $ B_{(j+1)\tau_k}$ the distance between
$x$ and any of "dangerous" subspaces will be greater than
$$
\frac{\rho_{(j+1)\tau_k}\gamma}{2} \ge \frac{\gamma}{4M^2t_{r_{j+2}}}\ge \frac{\gamma}{4M^{2+ k(\alpha,\beta, M)}t_{r_{j+1}} }.
$$
(in the last inequalities we use (\ref{t}) and (\ref{new})).

The inductive step is completed ant Lemma 1 is
proved.

 \vskip+1.0cm

author:

\vskip+0.3cm

Nikolay  G. Moshchevitin

e-mail: moshchevitin@mech.math.msu.su, moshchevitin@rambler.ru

\end{document}